\documentclass[notitlepage,leqno,11pt]{article}
\usepackage{amssymb,cite}
\catcode`\@=11 \@addtoreset{equation}{section}

\catcode`\@=12

\usepackage{latexsym}
\usepackage{textcomp}
\usepackage{amsmath}
\usepackage{amsfonts}
\usepackage{amssymb}
\usepackage{mathrsfs}

\renewcommand{\d}{\delta }

\newcommand{\D }{\Delta }

\newcommand{\e }{\varepsilon }
\newcommand{\g }{\gamma}

\renewcommand{\l }{\lambda }

\newcommand{\n }{\nabla }
\newcommand{\vp }{\varphi }
\renewcommand{\phi}{\varphi}

\newcommand{\s }{\sigma }

\renewcommand{\t }{\tau }

\renewcommand{\O }{\Omega }

\newcommand{\ov}{\overline}

\newcommand{\be}{\begin{equation}}
\newcommand{\ee}{\end{equation}}

\newcommand{\R}{\mathbb{R}}
\newcommand{\N}{\mathbb{N}}

\newcommand{\de}{\partial}

\newcommand{\ra}{{\rangle}}
\newcommand{\la}{{\langle}}

\newcommand{\calL }{\mathcal{L}}

\newcommand{\calD }{\mathcal{D}}

\newtheorem{Theorem}{Theorem}[section]
\newtheorem{Lemma}[Theorem]{Lemma}
\newtheorem{Proposition}[Theorem]{Proposition}
\newtheorem{Corollary}[Theorem]{Corollary}

\def\proof{\noindent{{\bf Proof. }}}
\def\square{\vbox{
    \hrule height .4pt
    \hbox{\vrule width .4pt height 7pt \kern 7pt
       \vrule width .4pt}
    \hrule height .4pt }}

\def\square{\vbox{
    \hrule height .4pt
    \hbox{\vrule width .4pt height 7pt \kern 7pt
       \vrule width .4pt}
    \hrule height .4pt }}

\def\QED{\hfill {$\square$}\goodbreak \medskip}

\def\R{{\mathbb R}}

\def\div{{\rm div}}
\def\eps{{\varepsilon}}

\newcommand{\Ds}{(-\Delta)^s}

\font\sc=cmcsc9  \textwidth=14truecm
\title{Nonexistence results for a class of fractional elliptic
  boundary value problems}
\author{Mouhamed Moustapha Fall   and   Tobias Weth}

\begin{document}
\date{}
\maketitle
 \let\thefootnote\relax\footnotetext{fall@math.uni-frankfurt (M. M. Fall), weth@math.uni-frankfurt (T. Weth).}
\let\thefootnote\relax\footnotetext{ Goethe-Universit\"{a}t Frankfurt, Institut f\"{u}r Mathematik.
Robert-Mayer-Str. 10 D-60054 Frankfurt, Germany.}
\bigskip

\begin{abstract}
In this paper we study a class of fractional elliptic problems of the form
$$
\left\{
\begin{aligned}
\Ds u&= f(x,u)&& \qquad \textrm{ in }\O\\
u&=0&&\qquad\textrm{ in }\R^N \setminus \O,
\end{aligned}
\right.
$$
where $s\in(0,1)$. We prove nonexistence of positive solutions when
$\O$ is star-shaped and $f$ is supercritical. We also derive a nonexistence result for
subcritical $f$ in some unbounded domains.
The argument relies on the method of moving spheres
applied to a reformulated problem using the
Caffarelli-Silvestre extension \cite{CSilv} of a solution of the above
problem. The standard approach in the case $s=1$ using Pohozaev type
identities does not carry over to the case $0<s<1$ due to the lack of
boundary regularity of solutions.

\end{abstract}

\section{Introduction}
Let $s\in (0,1)$ and $N>2s$. In the present paper, we are concerned
with the  nonexistence of positive functions
solving the fractional elliptic semilinear problem
\be
\label{eq:pblm}
\left\{
\begin{aligned}
\Ds u&=f(x,u)&&\quad\textrm{ in }\O,\\
u&=0  &&\quad\textrm{ in }\R^N\setminus\O.
\end{aligned}
\right.
\ee
in a domain $\Omega \subset \R^N$. Problems of this type received immensely growing attention recently,
while different versions of the nonlocal operator $\Ds$ related to
Dirichlet boundary conditions are
studied (see e.g. \cite{CT,CDDS,dPS,BCD,Tan,SV}). The version we consider in \eqref{eq:pblm} is the one most commonly
considered in analysis and probability theory. In probabilistic terms,
it can be defined as the generator of the $2s$-stable process in
$\Omega$ killed upon leaving $\Omega$. For our purposes, it is more
convenient to give an analytic definition. We define
$\Ds$ for any $\vp \in C^\infty_c(\R^N)$ by
\be\label{eq:DSvp-int}
\Ds \vp(x)= P.V.\int_{\R^N}\frac{\vp(x)-\vp(y)}{|x-y|^{N+2s}}dy= \lim_{\e\to0}\int_{|x-y|>\e}\frac{\vp(x)-\vp(y)}{|x-y|^{N+2s}}dy
\ee
for $x \in \R^N$, where P.V. stands for the principle value integral. We point out that this definition differs from the standard definition by a multiplicative constant. Via Fourier transform, \eqref{eq:DSvp-int}  is equivalent to
$$
C_{N,s} \widehat {\Ds \vp}(\xi)=|\xi|^{2s} \widehat \vp (\xi) \qquad \text{for
  $\xi \in \R^N$.}
$$
with the normalization constant $C_{N,s}=
s(1-s)\pi^{-N/2}2^{2s}\frac{\Gamma(\frac{N+2s}{2})}{\Gamma(2-s)}$, see
e.g. \cite[Remark 3.11]{CS}.
 Thanks to  Lemma \ref{lem:bondDsvp} below, for any
 $\vp \in C^\infty_c(\R^N)$ we have the  estimate
\begin{equation}
  \label{eq:1}
|\Ds\vp(x)|\le  C \frac{\|\vp\|_{C^2(\R^N)}}{1+|x|^{N+2s}}\quad \text{for all $x\in \R^N$,}
\end{equation}
where $C$ only depends on the support of $\vp$.
Let
$\calL^1_s$ denote the space of all measurable functions $u:\R^N\to\R$ such that
$\int_{\R^N}\frac{|u|}{1+|x|^{N+2s}}\,dx<\infty$, and let $\O$ be an open set of  $\R^N$.
 We define the Hilbert space $\calD^{s,2}(\O)$ as the completion of $C^\infty_c(\O)$ with respect to the norm
$\|\cdot\|_{\calD^{s,2}}$ induced by the scalar product $\langle \cdot ,\cdot  \rangle_{\calD^{s,2}}$ given by 
\be
\label{eq:norm}
\langle u,v \rangle_{\calD^{s,2}}= \int_{\R^{2N}}\frac{(u(x)-u(y))(v(x)-v(y))}{|x-y|^{N+2s}}dxdy.
\ee
We note that if $\Omega$ is a bounded Lipschitz domain, then $\calD^{s,2}(\O)$ coincides with the Sobolev space $\{u\in H^s(\R^N)\,:\, \textrm{$u=0$ a.e. in $\R^N\setminus\O$}\}$. We also observe that
 --  for any $u\in
 \calD^{s,2}(\O)$ -- the H\"{o}lder and the Hardy-Littlewood-Sobolev inequalities imply that
$$
\int_{\R^N}\frac{|u|}{1+|x|^{N+2s}} dx \leq C \|u\|_{\calD^{s,2}
} \qquad \text{for all $u \in \calD^{s,2}(\O)$}
$$
with a constant $C>0$. In other words, $\calD^{s,2}(\O)$ is
continuously embedded in $\calL^1_s$. As a
consequence, by recalling (\ref{eq:1}) we may define $\Ds u$ for every $u
\in \calD^{s,2}(\O)$ as a distribution by
$$
\la\Ds u,\vp \ra: =\int_{\R^N} u\Ds\vp\, dx=\langle u, \phi \rangle_{\calD^{s,2}} \qquad \text{for all $\vp \in
  C^\infty_c(\O).$}
$$
In particular, given $f\in L^1_{loc}(\O)$, we note that
$u\in\calD^{s,2}(\O)$ solves the problem $\Ds u=f$
if and only if
\begin{equation}
  \label{eq:2}
\langle u, \phi \rangle_{\calD^{s,2}} =\int_{\O}f\vp dx \qquad \text{for all $\vp\in C^\infty_c(\O).$}
\end{equation}
Throughout the paper, when we refer to solution of \eqref{eq:pblm}, we
mean distributional solutions $u\in\calD^{s,2}(\O)$ in the sense of (\ref{eq:2}) with
$f=f(\cdot,u(\cdot)) \in L^1_{loc}(\Omega)$. In order to state the main result of the
present paper, we need to introduce a definition of a star domain
which is slightly more general than usually considered in the
literature. We say that an open set  $\O \subset \R^N$ is star-shaped (or a
star domain) with respect to the origin $0 \in \overline \Omega$ if
for every $x \in \Omega$ we have $tx \in \Omega$ for $0<t \le 1$.
So in contrast to the standard definition, we also allow the star
center to lie on the boundary of $\Omega$. This will be crucial in
deriving results in unbounded domains. In particular, the punctured
open unit ball $B_1(0) \setminus \{0\}$ is star-shaped with respect to
the origin according to our definition. Our main result is the following.
\begin{Theorem}\label{th:neSupc}
Assume that $\O$ is bounded and star-shaped with respect to the origin
$0\in\ov{\O}$. Suppose that  $f: \ov{\Omega} \setminus \{0\} \times
[0,\infty) \to\R$
is locally Lipschitz in its second variable uniformly in compact subsets
of $ \ov{\Omega} \setminus\{0\}$  and
 is supercritical in the sense that
\be\label{eq:supCri}
\left\{
  \begin{aligned}
  &\text{the function $\l \mapsto \l^{-(N+2s)/(N-2s)} f(\l^{-2/(N-2s)}x,\l u)$}\\
&\text{is non-decreasing on $[1,\infty)$ for every $x \in \Omega
  \setminus \{0\},\,u \ge 0$}.
  \end{aligned}
\right.
\ee
Then \eqref{eq:pblm} has no positive solution $u\in
C(\R^N \setminus\{0\}) \cap \calD^{s,2}(\O)$.
\end{Theorem}

We remark that for $C^1$-nonlinearities $f: \ov{\Omega} \setminus \{0\} \times \R\to\R$
the supercriticality assumption \eqref{eq:supCri} is equivalent to
\begin{equation}
  \label{eq:3}
H_f(x,u) \ge 0\qquad \text{for all $(x,u) \in \Omega \setminus \{0\} \times [0,\infty)$,}
\end{equation}
where
\begin{equation}
  \label{eq:4}
H_f(x,u):= u \frac{\partial}{\partial u}f (x,u)- \frac{N+2s}{N-2s}f(x,u)-\frac{2}{N-2s}x
  \cdot \nabla_x f(x,u).
\end{equation}
As a first  consequence of Theorem \ref{th:neSupc}
we have the following Pohozaev type result.
\begin{Corollary}
\label{cor-poho}
Assume that $\O$ is bounded and star-shaped with respect to the
origin, and let $V\in C^1(\O \setminus \{0\})$ satisfy
$$
s V(x)+\frac{1}{2}\n V(x)\cdot x\geq 0\quad \text{for all $x\in \O
  \setminus \{0\}.$}
$$
Let $u\in \calD^{s,2}(\O)\cap C(\R^N \setminus \{0\})$, $u \geq 0$ in $\R^N$ be such that
\be\label{eq:potV}
\left\{
\begin{aligned}
 \Ds u+V(x)u&=u^p&&\qquad\textrm{ in }\O,\\
  u&=0 &&\qquad\textrm{ in } \R^N\setminus \O
\end{aligned}
\right.
\ee
for some $p\geq \frac{N+2s}{N-2s}$. Then $u = 0$ in $\R^N$.
\end{Corollary}

In the case where $\Omega$ is the unit ball in $\R^N$ and $V \equiv 0$, this gives an
affirmative answer to a conjecture of Birkner, L\'opez-Mimbela and
Wakolbinger, see \cite[p. 91]{BLW}. 
We note that existence results for problem \eqref{eq:potV} in the subcritical range $1<p< \frac{N+2s}{N-2s}$ and for more general
subcritical nonlinearities have been obtained recently by the first author
in \cite{Fall-frac} and by Servadei and Valdinoci in \cite{SV}.

In our next result the linear term is related to the relativistic Hardy inequality, see \cite{FLS} and \cite{Fall-frac}.
\begin{Corollary}
\label{cor:hardy}
 Assume that $\O$ is bounded and star-shaped with respect to the
 origin and let $u\in \calD^{s,2}(\O)\cap C(\R^N\setminus\{0\})$,
 $u\geq0$ in $\R^N$ be such that
\be\label{eq:pdHard}
\left\{
\begin{aligned}
\Ds u-\g |x|^{-2s}u&=u^p &&\quad\textrm{ in }\O,\\
u&=0&&\quad\textrm{ in }\R^N\setminus\O
\end{aligned}
\right.
\ee
for some $\g \in\R$ and $p\geq \frac{N+2s}{N-2s}$. Then $u=0$ in $\R^N$.
\end{Corollary}
Our next result is concerned
 with a singular nonlinearity.

\begin{Corollary}\label{cor:pdSW}
 Assume that $\O$ is bounded and star-shaped with respect to the
 origin and let $u\in \calD^{s,2}(\O)\cap C(\R^N\setminus\{0\})$,
 $u\geq0$ in $\R^N$ be such that
\be\label{eq:pdSW}
\left\{
\begin{aligned}
\Ds u&=|x|^{-\s}u^p&&\quad\textrm{ in }\O,\\
u&=0&&\quad\textrm{ in }\R^N\setminus \O
\end{aligned}
\right.
\ee
for some $\s \in \R$ and $p\geq \max \bigl\{1,\frac{N+2s-2\s}{N-2s}\bigr\}$. Then $u=0$ in $\R^N$.
\end{Corollary}

This result should be seen in the context of the criticality of $q=
\frac{2(N-\s)}{N-2s}=\frac{N+2s-2\s}{N-2s}+1$ for the embedding of the Sobolev space
$\calD^{s,2}(\O)$ in the weighted space
$L^{q}(\Omega;|x|^{-\s})$. More precisely, if $N>\max(\s,2s)$ and the underlying domain is
bounded, $\calD^{s,2}(\O)$ is continuously embedded in
$L^{q}(\Omega;|x|^{-\s})$ if and only if $q \le
\frac{2(N-\s)}{N-2s}$, and the embedding is compact iff
$q<\frac{2(N-\s)}{N-2s}$. Note also that the existence of the embeddings in
the subcritical range follows from the fact that
$$
\calD^{s,2}(\R^N) \hookrightarrow L^{2(N-\s)/(N-2s)}(\R^N;|x|^{-\s}),
$$
and this latter embedding can be seen as a version of the Stein-Weiss inequality \cite{StWe}.

Our next result is concerned with a class of unbounded
domains. Slightly extending a notion from \cite{RZ}, we say that an open set
 $\O$ is {\em star-shaped with
respect to infinity} if there exists a point $e\in\R^N\setminus\ov{\O}$
such that for every point $x\in \O$
the half-line $\{e+t(x-e)\,:\,t \ge 1\}$ is contained in
$\Omega$. Up to suitable translation, it is equivalent to require $0 \not \in
\overline \Omega$ and that $\R^N \setminus \overline \Omega$ is star-shaped
with respect to $0$ in the sense defined earlier.

\begin{Theorem}\label{th:ne-stInf}
Assume that $\O$ is  star-shaped with respect to infinity.
Let $u\in \calD^{s,2}(\O)\cap C(\R^{N})$ be nonnegative and such that
\be\label{eq:ustf-stU}
\left\{
\begin{aligned}
\Ds u&= u^p &&\qquad\textrm{in }\O,\\
u&=0 &&\qquad\textrm{in }\R^N\setminus\O
\end{aligned}
\right.
\ee
for some $1\le  p\leq \frac{N+2s}{N-2s}$. Then $u = 0$ in $\R^N$.
\end{Theorem}

In fact, we will deduce Theorem~\ref{th:ne-stInf} from
Theorem~\ref{th:neSupc} via a variant of the classical Kelvin
transform, see Sections~\ref{sec:some-preliminaries}
and~\ref{sec:proof-main-results} below for
details.\\
 Theorem~\ref{th:ne-stInf} in particular applies to the
cone-like domains $\Omega_\tau:= \{x \in \R^N \setminus \{0\}\::\:
\frac{x_N}{|x|}>\tau \}$ for $\tau \in (-1,1)$. Here one may take
$e=-e_N$, where $e_N$ is the $N$-th coordinate vector, in the
definition of star-shapedness at infinity. Since the half-space
$\R^N_+$ is a particular case with $\tau=0$, we deduce the following corollary.

\begin{Corollary}\label{cor:ne-hs}
Let $u\in \calD^{s,2}(\R^{N}_+)\cap C(\R^{N})$ be nonnegative and such
that
\be\label{eq:ustf-hs}
\left\{
\begin{aligned}
 \Ds u&= u^p &&\qquad\textrm{in }\R^{N}_+,\\
u&=0 &&\qquad\textrm{in }\R^N\setminus\R^{N}_+
\end{aligned}
\right.
\ee
for some $1\leq p\leq \frac{N+2s}{N-2s}$. Then $u= 0$ in $\R^N$.
\end{Corollary}

We remark that Theorem~\ref{th:ne-stInf} does not apply to the case $\O=\R^N$.
Indeed, in this case the critical problem with $p=\frac{N+2s}{N-2s}$ admits positive solutions which have been classified completely in \cite{CLO}. Moreover, in the case $\Omega= \R^N$, $s\in[1/2,1)$ and $1< p<\frac{N+2s}{N-2s}$, a nonexistence result has been obtained very recently and independently in \cite{dPS} by de Pablo and S\'anchez, see also \cite{H} for $s=1/2$ and  $1< p<\frac{N+2s}{N-2s}$.  

In order to explain our approach to obtain the nonexistence results, we
need to compare \eqref{eq:pblm} with the classical problem
\be
\label{eq:pblm-s=1}
\left\{
\begin{aligned}
-\Delta u&=f(x,u)&&\qquad\textrm{ in }\O,\\
u&=0  &&\qquad\textrm{on }\partial \O.
\end{aligned}
\right.
\ee
For \eqref{eq:pblm-s=1}, the analogue of Theorem~\ref{th:neSupc} is true, and for strictly starshaped $C^1$-domains
$\Omega$ and
$C^1$-nonlinearities $f$ on $\Omega \times [0,\infty)$ satisfying
additionally $f(\cdot,0)=0$ it can be derived from the Pohozaev type integral identity
\begin{equation}
  \label{eq:10}
\int_\Omega \int_0^{u(x)} H_f(x,t)\,dt \,dx +\frac{1}{N-2}
\int_{\partial \Omega} u_\nu^2\: x \!\cdot\! \nu\,d \sigma(x) = 0,
\end{equation}
see e.g. \cite[Theorem 5.2]{Reichel}. Here $H_f$ is defined as in
(\ref{eq:4}). Indeed, by (\ref{eq:3}) and the star-shapedness of
$\Omega$, the LHS of \eqref{eq:10} is nonnegative, and  by unique continuation
it is strictly positive if $u \not \equiv 0$. The above integral
identity can be derived by multiplying \eqref{eq:pblm-s=1} with the functions $u$ and $x \mapsto x \cdot \nabla u$ respectively and integrating
by parts. The same strategy does not work for
(\ref{eq:pblm}) since the problem is nonlocal and does not allow a
simple integration by parts formula as in the case $s=1$. More severely,
in the case $0<s<1$ solutions of~\eqref{eq:pblm} are not of class
$C^1$ up to the boundary  even if the underlying domain is smooth. In particular, if $x \mapsto f(x,u(x)) \ge 0$ is a nonnegative nontrivial
function on $\Omega$, then any solution $u$ of \eqref{eq:pblm} fails to
possess a finite normal derivative $u_\nu$ on $\partial \Omega$, see
e.g. \cite[Lemma 4.3]{BLW}.\\
The approach we follow here is inspired
by Reichel and Zou \cite{RZ} who used the technique of moving spheres
to prove nonexistence results for cooperative elliptic systems. The
moving sphere method can be seen as a variant of the method of moving
hyperplanes (see e.g. \cite{A,S,GNN1,GNN2,BN}) and has been widely
used to classify positive solutions of nonlinear elliptic problems, see
e.g. \cite{YanYanLi} and the references therein. For the special case where the underlying domain is the entire space $\R^N$, it has also been applied to problems involving the fractional Laplacian, 
see the aforementioned recent paper \cite{dPS} of de Pablo and S\'anchez and also \cite{CLO}. Unlike as in \cite{RZ}, we are not able to implement a moving sphere argument directly in the present setting, so instead -- as in \cite{dPS} --  we first transform \eqref{eq:pblm} to a local problem by
considering the Caffarelli-Silvestre extension of a solution $u$
on $\R^{N+1}_+$, see \cite{CSilv} and also \cite{CS,Fall-frac}.
This extension satisfies $w =u$
on $\O$ and solves in some weak sense (see
Section~\ref{sec:some-preliminaries} for details) the boundary value problem
\be\label{eq:Harmext}
\left\{
\begin{aligned}
\div(t^{1-2s}\n w)&=0 &&\quad \textrm{ in }\R^{N+1}_+,\\
w& =0&& \quad \textrm{ on }\R^N\setminus\O,\\
- c_{N,s} \lim \limits_{t \to 0^+}t^{1-2s}\,\frac{\de w}{\de t}&= f(x,w) &&  \quad \textrm{ on } \O,
\end{aligned}
\right.
\ee
with the positive normalization constant $c_{N,s}=\frac{\pi^{N/2} \Gamma(s)}{2s
  \Gamma(\frac{N+2s}{2})}$ (note that this constant is different from
the one noted e.g. in \cite[Remark 3.11]{CS} due to our
normalization of $(-\Delta)^s$). Here and in the following we write $z=(x,t) \in \R^{N+1}_+$ with
$x \in \R^N$ and $t>0$, and we identify $\R^N$ with $\partial \R^{N+1}_+$, so that
$\Omega$ is contained in $\partial \R^{N+1}_+$.  We will then
apply the moving sphere argument to the local problem
\eqref{eq:Harmext} in place of \eqref{eq:pblm}. We note that the
Caffarelli-Silvestre extension of a solution of \eqref{eq:pblm} has
received considerable attention in recent years due to its usefulness
in the context of many different problems, see e.g. \cite{CSS,SirV,DS,CG,CRS}.

We should mention  that  -- in contrast to the
nonexistence results for \eqref{eq:pblm-s=1} based on the Pohozaev type
identity -- our approach does not extend to sign changing solutions.
The existence resp. nonexistence of sign changing solutions of
\eqref{eq:pblm} under supercriticality and star-shapedness assumptions
therefore remains an open problem.

Finally, we would like to compare \eqref{eq:pblm} with the related problem
\be
\label{eq:pblm-spectral}
\left\{
\begin{aligned}
A^s u&=f(x,u)&&\quad\textrm{ in }\O,\\
u&=0  &&\quad\textrm{on $\partial \O$.}
\end{aligned}
\right.
\ee
Here $A$ stands for the negative Laplacian as a self adjoined operator in
$L^2(\Omega)$ with domain
$$
\{u \in H^1_0(\Omega)\::\: \text{$\Delta u \in L^2(\Omega)$ as a
  distribution}\},
$$
and $A^s$ is the corresponding power in spectral theoretic
sense. Although problems \eqref{eq:pblm} and \eqref{eq:pblm-spectral}
look similar, there are crucial differences as discussed e.g. in \cite{Fall-frac}. In particular, solutions
of \eqref{eq:pblm-spectral} have in general much better
boundary regularity than solutions of \eqref{eq:pblm}, and this can
also be seen when comparing the corresponding extended problems. We
point out that in
\cite{CT,CDDS,dPS,BCD,Tan} a variant of the Caffarelli-Silvestre
extension for solutions of \eqref{eq:pblm-spectral} was  considered
which preserves the regularity properties up to the
boundary. Moreover, nonexistence results for \eqref{eq:pblm-spectral}
have recently been
proved in \cite{Tan,BCD} via a Pohozaev
type integral identity for the extended problem. As we pointed out
before, such an approach is not available for \eqref{eq:pblm}
resp. \eqref{eq:Harmext} due to the lack of boundary regularity of solutions.

The paper is organized as follows. In
Section~\ref{sec:some-preliminaries} we discuss a suitable weak notion
of solution of \eqref{eq:Harmext}, and we study how
problems~\eqref{eq:pblm} and \eqref{eq:Harmext} transform under a Kelvin
type transform. We also formulate two versions of boundary
maximum principles related to a linearized  version of problem
\eqref{eq:Harmext}. Since this section deals with all technical
aspects of the problem, the remaining parts of the proofs of our main
results are
relatively short, and they are contained in Section~\ref{sec:proof-main-results}.

\medskip

\noindent
\textbf{Acknowledgments:}
This work is supported by the Alexander von Humboldt foundation. The
authors would like to thank the referee for his/her valuable remarks.

\section{Some preliminaries}
\label{sec:some-preliminaries}

Throughout the paper, we consider $s\in(0,1)$ and assume that $N>2s$. In this section we collect preliminary tools related to
\eqref{eq:pblm} and the reformulated version \eqref{eq:Harmext}. We
also need to introduce some definitions concerning notions
of weak solutions. If $\Omega \subset \R^N$ is an open set and $f \in
L^1_{loc}(\Omega)$, we say that $u \in \calD^{s,2}(\R^N)$ is a
distributional solution of $(-\Delta)^s u = f$ in $\Omega$ if
\begin{equation}
  \label{eq:7}
\langle u,\vp \rangle_{\calD^{s,2}} = \int_{\Omega} f \vp \,dx \qquad \text{for all $\vp \in C^\infty_c(\Omega),$}
\end{equation}
where $\langle \cdot,\cdot \rangle_{\calD^{s,2}}$ is defined in
(\ref{eq:norm}).  Note that by considering $u \in \calD^{s,2}(\R^N)$
we do not prescribe $u$ on $\R^N \setminus \Omega$ here. We start with the following result.
\begin{Lemma}\label{lem:bondDsvp}
Let $\O$ be a bounded open set.  Then there
 exists a constant $C=C(N,s,\O)>0$ such that for all $\vp\in
 C^2_c(\O),\: x\in \R^N$ and $\e\in(0,1)$ we have

 \begin{equation}
   \label{eq:24}
\left|\int_{|x-y|>\e}\frac{\vp(x)-\vp(y)}{|x-y|^{N+2s}}dy  \right|\leq
\frac{C\|\vp\|_{C^2(\R^N)}}{1+|x|^{N+2s}}  \end{equation}

\end{Lemma}
\proof
For $x \in \R^N$ and $\eps>0$, integration by parts yields
\begin{align*}
\int_{|x-y|>\e}&\frac{\vp(x)-\vp(y)}{|x-y|^{N+2s}}dy=\int_0^1
\int_{|x-y|>\e}\n\vp (x+t(y-x))\cdot \frac{x-y}{|x-y|^{N+2s}}dy dt\\
&=\frac{1}{N+2s-2} \Bigl(\int_0^1 \int_{|x-y|=\e}
\n\vp(x+t(y-x))\cdot(y-x)|x-y|^{-N-2s+1}d\s(y) dt\\
&+\int_0^1 t \int_{|x-y|>\e} \D\vp(x+t(y-x))|x-y|^{-N-2s+2}dydt\Bigr),
\end{align*}
whereas
\begin{align*}
\int_0^1& \int_{|x-y|=\e} \n\vp(x+t(y-x))\cdot(y-x)|x-y|^{-N-2s+1}d\s(y) dt\\
&=\e^{1-2s}\int_0^1\int_{S^{N-1}}\n\vp(x+t\e\s)\cdot\s d\s dt\\
&=\e^{1-2s}\int_0^1\int_{S^{N-1}}\n\vp(x)\cdot\s d\s dt
+ \e^{2-2s}\int_0^1t\int_0^1 \int_{S^{N-1}} D^2\vp(x+\e
t\t\s)[\s]\cdot\s d\s d\t dt
\end{align*}
and, by oddness,
$$
\int_{S^{N-1}}\n\vp(x)\cdot\s d\s=\sum_{i=1}^N\frac{\de\vp}{\de x^i}(x) \int_{S^{N-1}}\s^i d\s =0.
$$
Consequently,
\begin{align}
\int_{|x-y|>\e}\frac{\vp(x)-\vp(y)}{|x-y|^{N+2s}}dy &=
\frac{1}{N+2s-2} \int_0^1 t \int_{|x-y|>\e} \D\vp(x+t(y-x))|x-y|^{-N-2s+2}dydt
\nonumber\\
&+\frac{\e^{2(1-s)}}{N+2s-2}\int_0^1t\int_0^1 \int_{S^{N-1}}
D^2\vp(x+\e t\t\s)[\s]\cdot\s d\s d\t dt, \label{eq:22}
\end{align}
while
\begin{equation}
  \label{eq:20}
\Bigl|\int_0^1t\int_0^1 \int_{S^{N-1}}
D^2\vp(x+\e t\t\s)[\s]\cdot\s d\s d\t dt\Bigr| \le C_1 \|\phi\|_{C^2(\R^N)}
\end{equation}
with a constant $C_1>0$ depending only on $N$ and $s$. We now fix
$R>0$ such that $\O\subset B(0,R)$, and we first consider $x \in \R^N\setminus
  B(0,4R)$. Then $|x-y|\geq R+\frac{|x|}{2}$ for $y \in \Omega$ and therefore
\begin{equation}
  \label{eq:23}
\left|\int_{|x-y|>\e}\frac{\vp(x)-\vp(y)}{|x-y|^{N+2s}}dy\right|\leq
 \int_{|y| \le R} \frac{|\vp(y)|}{(R+\frac{|x|}{2})^{N+2s}}dy
\leq C_2 \frac{\|\vp\|_{C^2(\R^N)}}{1+|x|^{N+2s}}
\end{equation}
with a constant $C_2>0$ depending only on $R$, $N$ and $s$. Next we
consider $x \in B(0,4R)$ and note that, for every $t \in (0,1)$,
$$
|x-y|\leq  \frac{R+|x|}{t}\le \frac{5R}{t} \qquad \text{if $|x+t(y-x) |\leq R$},
$$
and
$$
\D\vp(x+t(y-x)) = 0 \qquad \text{if $|x+t(y-x)| \ge R$.}
$$
Hence for $x \in B(0,4R)$ we have
\begin{align}
\Big| \int_0^1 \!t \! \int_{|x-y|>\e}
\!\!\!\!\!\D\vp(x+& t(y-x))|x-y|^{-N-2s+2}dydt \Big|
\nonumber\\
& \le
\|\vp\|_{C^2(\R^N)}\!\int_0^1\!\!t\!\int_{|x+t(y-x)|<R}\!\!\!|x-y|^{-N-2s+2}dydt
\nonumber\\
& \le
\|\vp\|_{C^2(\R^N)}\!\int_0^1\!\!t\!\int_{|x-y|\le \frac{5R}{t}}|x-y|^{-N-2s+2}dydt
\nonumber\\
&\le \|\vp\|_{C^2(\R^N)}|S^{N-1}|\!\int_0^1\!\!t\!\int_0^{\frac{5R}{t}}\!\!\!r^{1-2s}drdt
\nonumber\\
&= \|\vp\|_{C^2(\R^N)} |S^{N-1}|\frac{(5R)^{2s-2}}{2-2s} \int_0^1 t^{-1+2s}dt= C_3
\|\vp\|_{C^2(\R^N)} \label{eq:21},
\end{align}
with a constant $C_3>0$ depending only on $R$, $N$ and $s$. Combining
(\ref{eq:22}), (\ref{eq:20}), (\ref{eq:23}) and (\ref{eq:21}), we find that there
exists a constant $C>0$ depending only on $R'$, $N$ and $s$ such that
(\ref{eq:24}) holds, as claimed.
\QED

Next, we consider the conformal diffeomorphism
\begin{equation}
  \label{eq:15}
\kappa: \R^N \setminus \{0\} \to \R^N \setminus \{0\},\qquad \kappa(x)=\frac{x}{|x|^2}.
\end{equation}
It is easy to see that
\begin{equation}
  \label{eq:16}
|\kappa(x)-\kappa(y)|=\frac{|x-y|}{|x||y|} \qquad \text{for every $x,y \in \R^N \setminus \{0\}$},
\end{equation}
and that the Jacobian determinant of $\kappa$ satisfies
$$
|\det J_\kappa(x)|=|x|^{-2N}.
$$
In the following, for a measurable function $u$ on $\R^N$, we
a.e. define $Ku$ on $\R^N$ by $$Ku(x)=|x|^{2s-N}u\left(\kappa(x) \right).$$
The map $K$ is usually called {\em Kelvin transform}, and it is a well known tool in potential theory and partial differential equations. It has also been studied in detail in a probabilistic framework for stable processes, 
see \cite{BZ} and the references therein. Here we need the following property of $K$.
  
\begin{Lemma}\label{eq:chang}
The map $K$ defines an isometry on $\calD^{s,2}(\R^N)$, i.e. for every
$u,v \in \calD^{s,2}(\R^N)$ we have $Ku,Kv \in \calD^{s,2}(\R^N)$ and
\begin{equation}
  \label{eq:5}
\langle u, v \rangle_{\calD^{s,2}} = \langle Ku, Kv \rangle_{\calD^{s,2}}.
\end{equation}
 \end{Lemma}
\proof
Since $C^\infty_c(\R^N\setminus\{0\})$ is dense in
$C^\infty_c(\R^N)$ with respect to the $\calD^{s,2}(\R^N)$-norm
as a consequence of our general assumption $N>2s$ (see \cite[p. 397]{Maz}), it suffices to show (\ref{eq:5}) for
$u,v\in C^\infty_c(\R^N\setminus\{0\})$. By changing variables and using~(\ref{eq:16}), we have
$$
\begin{array}{c}
\langle u, v \rangle_{\calD^{s,2}} =\displaystyle \int_{\R^{2N}}\frac{(u(x)-u(y))(v(x)-v(y))}{|x-y|^{N+2s}}dxdy\hspace{7cm} \vspace{3mm}\\
\displaystyle =\int_{\R^{2N}}\frac{(u(\kappa(x))-u(\kappa(y)))(v(\kappa(x))-v(\kappa(y)))}{|x-y|^{N+2s} |x|^{-N-2s}|y|^{-N-2s}}|x|^{-2N}|y|^{-2N} dxdy\hspace{0.9cm}\vspace{3mm}\\
\displaystyle= \int_{\R^{2N}}\frac{(u(\kappa(x))-u(\kappa(y)))(v(\kappa(x))-v(\kappa(y)))}{|x-y|^{N+2s} }|x|^{-N+2s}|y|^{-N+2s} dxdy.
\end{array}
$$
Observe that
$$
\begin{array}{c}
\displaystyle(u(\kappa(x))-u(\kappa(y)))(v(\kappa(x))-v(\kappa(y)))|x|^{-N+2s}|y|^{-N+2s}\hspace{5cm}\vspace{3mm}\\
\displaystyle= (Ku(x)-Ku(y))(Kv(x)-Kv(y))+ Ku(x)v(\kappa(x))[|y|^{2s-N}-|x|^{2s-N} ]\vspace{3mm}\\
\displaystyle + Ku(y)v(\kappa(y))[|x|^{2s-N}-|y|^{2s-N} ].
\end{array}
$$
We therefore have
$$
\begin{array}{c}
\displaystyle\int_{\R^{2N}}\frac{(u(x)-u(y))(v(x)-v(y))}{|x-y|^{N+2s}}dxdy
\displaystyle=\int_{\R^{2N} }\frac{  (Ku(x)-Ku(y))(Kv(x)-Kv(y)) }{ |x-y|^{N+2s}  }dxdy\vspace{3mm}\\
\displaystyle\hspace{3.2cm}+2 \lim_{\e\to0}\int_{\R^N}\int_{|x-y|>\e}\frac{Ku(x)v(\kappa(x))[|y|^{2s-N}-|x|^{2s-N} ]}{|x-y|^{N+2s}}dydx.
\end{array}
$$
It thus remains to prove that
\begin{equation}
  \label{eq:18}
\lim_{\e\to0}\int_{\R^N}\int_{|x-y|>\e}\frac{Ku(x)v(\kappa(x))[|y|^{2s-N}-|x|^{2s-N} ]}{|x-y|^{N+2s}}dydx=0.
\end{equation}
To show this, we consider $f\in C^\infty_c(\R^N\setminus\{0\})$
defined by $f(x)=Ku(x)v(\kappa(x))$. Since
$$\int_{\R^N}\int_{|x-y|>\e}\frac{f(x)|y|^{2s-N}}{|x-y|^{N+2s}} dydx <\infty,\quad \int_{\R^N}\int_{|x-y|>\e}\frac{f(x)|x|^{2s-N} }{|x-y|^{N+2s}} dydx  <\infty,$$
 we have by Fubini's theorem
$$
\int_{\R^N}\int_{|x-y|>\e}\frac{f(x)[|y|^{2s-N}-|x|^{2s-N} ]}{|x-y|^{N+2s}} dydx=\int_{\R^N}\int_{|x-y|>\e}\frac{|x|^{2s-N}(f(y)-f(x))}  {|x-y|^{N+2s}}dydx.
$$
Note that $x\mapsto |x|^{2s-N}\in \calL^1_s $. By Lemma \ref{lem:bondDsvp}, we have
$$
\left|\int_{|x-y|>\e}\frac{f(x)-f(y)}  {|x-y|^{N+2s}}dy\right|\leq
\frac{C}{1+|x|^{N+2s}}\quad \text{for all $\e\in(0,1)$},
$$
and therefore the dominated convergence theorem implies that
$$
\lim_{\e\to 0}\int_{\R^N}\int_{|x-y|>\e}\frac{f(x)[|y|^{2s-N}-|x|^{2s-N} ]}{|x-y|^{N+2s}} dydx=\int_{\R^N} |x|^{2s-N} \Ds f(x)dx.
$$
Since $x \mapsto |x|^{2s-N}$ is the Riesz potential of order $2s$, we have  (up to a constant)
$$
\int_{\R^N} |x|^{2s-N} \Ds f(x)dx=\la\Ds|x|^{2s-N},f \ra=\la \d, f \ra=0
$$
in distributional sense, because $f$ is supported away from the origin and $\d$ is the Dirac mass at the origin. Hence we have proved (\ref{eq:18}) and the lemma then follows.
\QED

As a consequence, we get the following result, which is closely related to \cite[Theorem 2]{BZ}. We note that, unlike in the present paper, probabilistic techniques are used in \cite{BZ}.

\begin{Corollary}
\label{sec:some-preliminaries-2}
Let $\Omega \subset \R^N$ be an open set and
$$
\tilde \Omega:=\kappa(\Omega
\setminus \{0\}) \subset \R^N \setminus \{0\}.
$$
Let $f \in
L^1_{loc}(\Omega)$, and let $u \in
\calD^{s,2}(\R^N)$ solve $(-\Delta)^s u= f$ in $\Omega$ in
distributional sense. Then $\tilde u = K u$ is contained in
$\calD^{s,2}(\R^N)$ and solves $(-\Delta)^s \tilde u = \tilde f$ in
distributional sense in $\tilde \Omega$, where $\tilde f \in
L^1_{loc}(\tilde \Omega)$ is given by $\tilde
f(x)=|x|^{-(N+2s)}f(\frac{x}{|x|^2})$.\\
Moreover, if $u \in \calD^{s,2}(\Omega)$, then $\tilde u \in
\calD^{s,2}(\tilde \Omega)$.
\end{Corollary}

\proof
Suppose first that $u \in \calD^{s,2}(\Omega)$. Since, as noted before, $C^\infty_c(\Omega
\setminus \{0\})$ is dense in
$\calD^{s,2}(\Omega)$, there exists a sequence $(\psi_n)_n$ in $C^\infty_c(\Omega
\setminus \{0\})$ with $\|u-\psi_n\|_{\calD^{s,2}} \to 0$ as $n \to
\infty$. By \eqref{eq:chang}, we then also have $\|Ku-K\psi_n\|_{\calD^{s,2}} \to 0$ as $n \to
\infty$. Since $K \psi_n \in \calD^{s,2}(\tilde \Omega)$ for all
$n$, this implies $K u \in \calD^{s,2}(\tilde \Omega)$.\\
Next we assume that $u \in \calD^{s,2}(\R^N)$ solves $(-\Delta)^s u= f$ in $\Omega$ in
distributional sense. Applying the
argument above to $\Omega= \R^N$ yields $\tilde u \in \calD^{s,2}(\R^N
\setminus \{0\}) \subset \calD^{s,2}(\R^N)$. Moreover, for given
$\tilde \phi
\in C^\infty_c(\tilde \Omega)$, we may now write $\tilde \phi= K
\phi$ with $\phi \in C^\infty_c(\tilde \Omega)$. By Lemma~\ref{eq:chang}, we
then have
\begin{align*}
\langle \tilde u, \tilde \phi \rangle_{\calD^{s,2}} =
\langle u, \phi \rangle_{\calD^{s,2}} &= \int_{\Omega} f \phi\,dx= \int_{\tilde \Omega} (f \circ \kappa ) (\phi \circ
\kappa) |\det J_\kappa|\,dx\\
&=\int_{\tilde \Omega}f(\kappa(x))
\phi(\kappa(x))|x|^{-2N}  \,dx= \int_{\tilde \Omega}\tilde f  \tilde \phi \,dx.
\end{align*}
This shows the claim.
\QED

Next, we introduce some notations related to the reformulated
problem~\eqref{eq:Harmext}. As before, we write $z=(x,t) \in
\R^{N+1}_+$ with $x \in \R^N$ and $t \in (0,\infty)$. Let
$D^{1,2}(\R^{N+1}_+;t^{1-2s})$ denote the space
of all functions $w \in H^1_{loc}(\R^{N+1}_+)$ such that
$$
\int_{\R^{N+1}_+} t^{1-2s} |\nabla w|^2\,dz<\infty.
$$
Formally introducing the operator $L_s:= \div (t^{1-2s} \nabla)$ on
$\R^{N+1}_+$, we say that a function $w \in D^{1,2}(\R^{N+1}_+;t^{1-2s})$ is
{\em weakly $L_s$-harmonic} if
$$
\int_{\R^{N+1}_+} t^{1-2s} \nabla w \nabla \phi\,dz=0 \qquad
\text{for all $\phi \in C^\infty_c(\R^{N+1}_+)$.}
$$
By standard elliptic regularity, every weakly $L_s$-harmonic function
$w \in D^{1,2}(\R^{N+1}_+;t^{1-2s})$ belongs to
$C^\infty(\R^{N+1}_+)$ and satisfies $\div (t^{1-2s} \nabla w) \equiv
0$ pointwise in $\R^{N+1}_+$. Moreover, $w$ does not attain an
interior maximum or minimum point in $\R^{N+1}_+$ unless $w$ is
constant. Note also that we have a well defined continuous  trace map
$$
D^{1,2}(\R^{N+1}_+;t^{1-2s}) \to \calD^{s,2}(\R^N)
$$
(see e.g. \cite{BCD}), and for the
sake of simplicity we denote the trace of a function in
$D^{1,2}(\R^{N+1}_+;t^{1-2s})$ with the same letter as the function
itself. If
$\phi,\psi \in D^{1,2}(\R^{N+1}_+;t^{1-2s})$ and $\phi$ is weakly
$L_s$-harmonic, we have the identity
\begin{equation}
  \label{eq:11}
c_{N,s} \int_{\R^{N+1}_+} t^{1-2s} \nabla \phi \nabla \psi \,dz=
\int_{\R^{2N}}\frac{
(\phi(x)-\phi(y))(\psi(x)-\psi(y))}{|x-y|^{N+2s}}dxdy
\end{equation}
with $c_{N,s}$ as in (\ref{eq:Harmext}). Now, for an open set $\Omega \subset \R^N$, we denote by
$D(\Omega,s)$ the closed subspace of functions in $D^{1,2}(\R^{N+1}_+;t^{1-2s})$ such that their trace on $\R^N$ is contained in $\calD^{s,2}(\O)$.
 It is easy to see that every function
$u \in \calD^{s,2}(\Omega)$
has a unique weakly harmonic extension $H(u) \in
D(\Omega,s)$ which can be found by minimizing the functional
$$
w \mapsto \int_{\R^{N+1}_+} t^{1-2s} |\nabla w|^2\,dz
$$
among all functions $w \in D(\Omega,s)$ satisfying
$w=u$ on $\R^N$. Using this fact in the special case $\Omega= \R^N$
(in which $D(\R^N\!,s)=D^{1,2}(\R^{N+1}_+;t^{1-2s})$) together with (\ref{eq:11}), we find that
\begin{equation}
  \label{eq:12}
\int_{\R^{2N}}\frac{
(\phi(x)-\phi(y))^2}{|x-y|^{N+2s}}dxdy \le c_{N,s} \int_{\R^{N+1}_+} t^{1-2s}
|\nabla \phi|^2\,dz
\end{equation}
for all $\phi \in D^{1,2}(\R^{N+1}_+;t^{1-2s})$. Moreover, since $\calD^{s,2}(\R^N)$ is continuously embedded in
$L^{\frac{2N}{N-2s}}(\R^N)$, there exists a constant $C>0$ such that
\begin{equation}
  \label{eq:13}
\|\phi\|_{L^{\frac{2N}{N-2s}}(\R^N)}^2 \le C \int_{\R^{2N}}\frac{(\phi(x)-\phi(y))^2}{|x-y|^{N+2s}}dxdy \qquad \text{for all $\phi \in \calD^{s,2}(\R^N)$.}
\end{equation}
Another fact we need
is the following:

\begin{Lemma}
\label{sec:some-preliminaries-3}
Let $\Omega \subset \R^N$ be a bounded open set, and let $u \in
D^{1,2}(\R^{N+1}_+;t^{1-2s})$ be such that its trace -- also denoted
by $u$ -- is continuous
in $\overline \Omega$ and satisfies $u \equiv 0$ on $\R^N \setminus \Omega$. Then $u \in D(\Omega,s)$.
\end{Lemma}

\proof
Consider $G\in C^\infty(\R)$  such that
$$
 G(r)=0\quad \textrm{ if } |r|\leq 1,\quad G(r)=r\quad \textrm{
   if }|r|\geq 2\quad\textrm{ and }\quad |G'(r)|\leq 3 \quad
 \text{if $1 \le |r| \le 2$.}
$$
Then the functions $u_n$ defined by $u_n(t,x)=\frac{1}{n}G(nu(t,x))$
 are clearly contained in $D^{1,2}(\R^{N+1}_+;t^{1-2s})$ for $n \in
 \N$. Passing to traces, we therefore have
$u_n\in \calD^{s,2}(\R^N)$. Note that by the dominated convergence
theorem we have $u_n\to u $ in $ D^{1,2}(\R^{N+1}_+;t^{1-2s})$.
In addition, since the support of the trace of $u_n$ in $\R^N$, is contained in the compact subset of $\O$
$$
\left\{x\in\O\,:\, |u_n(x)|\geq \frac{1}{n}\right\},
$$
it follows that $u_n\in  \calD^{s,2}(\O)$ by the density result in \cite[Theorem 1.4.2.2]{Gris}. To conclude we observe that $u_n\to u$ in $  \calD^{s,2}(\O)$
and this holds true thanks to the continuity of the trace operator  $ D^{1,2}(\R^{N+1}_+;t^{1-2s})\to \calD^{s,2}(\R^N)$.
\QED
We remark that the continuity assumption in {Lemma} \ref{sec:some-preliminaries-3}   is not needed if $\O$ has a continuous boundary, see \cite[Theorem 1.4.2.2]{Gris}.\\

 Next, let $q_s:=
\frac{2N}{N+2s}$ be the conjugate of $\frac{2N}{N-2s}$. If $f \in L^{q_s}(\Omega)$
is given and $u \in \calD^{s,2}(\R^N)$ satisfies $(-\Delta)^s u=
f$ in $\Omega$ in distributional sense, then, as a consequence of the
embedding $\calD^{s,2}(\Omega) \hookrightarrow
L^{\frac{2N}{N-2s}}(\Omega)$, it also satisfies
this equation in weak sense, i.e.
$$
\langle u, \psi \rangle_{\calD^{s,2}}= \int_{\Omega} f
\psi\,dx
\qquad \text{for all $\psi \in \calD^{s,2}(\Omega)$.}
$$
Moreover, by (\ref{eq:11}), the weakly $L_s$-harmonic extension
$w=H(u) \in D(\Omega,s)$ of $u$ then satisfies
\begin{equation}
  \label{eq:14}
c_{N,s} \int_{\R^{N+1}_+} t^{1-2s} \nabla w \nabla \psi \,dz =
\int_{\Omega}f \psi\,dx \qquad \text{for all $\psi \in D(\Omega,s)$.}
\end{equation}
We may summarize the discussion in the following statement.

\begin{Lemma}
\label{sec:some-preliminaries-4}
Let $\Omega \subset \R^N$ be an open set and $f \in
L^{q_s}(\Omega)$. A function $w \in
D^{1,2}(\R^{N+1}_+;t^{1-2s})$ satisfies (\ref{eq:14}) if and only if
$w$ is weakly $L_s$-harmonic and its trace -- also denoted by $w \in
\calD^{s,2}(\R^N)$ -- solves $(-\Delta)^s w = f$ in $\Omega$ in
distributional sense.\\
If this holds, we say that $w$ weakly solves the problem
\begin{equation}
\label{eq:19}
\left \{
\begin{aligned}
\div (t^{1-2s} \nabla w)&=0&&\quad \text{in $\R^{N+1}_+$,}\\
-c_{N,s}\lim \limits_{t \to 0}t^{1-2s}w_t&=f&&\quad \text{on $\Omega$.}
\end{aligned}
\right.
\end{equation}
\end{Lemma}

Next, we examine how problems of type (\ref{eq:19}) transform under
generalized Kelvin inversions.

\begin{Proposition}
\label{sec:some-preliminaries-1}
Let $w \in D^{1,2}(\R^{N+1}_+;t^{1-2s})$, let $\Omega \subset
\R^N$ be an open set and let $f \in
L^{q_s}(\Omega)$. Moreover, for fixed $\rho>0$, consider
$$
\Omega_\rho:= \Bigl \{\frac{\rho^2x}{|x|^2}\::\: x \in \Omega \setminus
\{0\}\Bigr \} \subset \R^N,
$$
and let $w_\rho: \R^{N+1}_+ \to \R$, $f_\rho: \Omega_\rho \to \R$ be defined by
$$
w_\rho(z):=\left(\frac{\rho}{|z|}\right)^{N-2s}w\left(\frac{\rho^2
    z}{|z|^2}\right) \qquad \text{and}\qquad
f_\rho(x)= \left(\frac{\rho}{|x|}\right)^{N+2s}f\left(\frac{\rho^2 x}{|x|^2}\right).
$$
Then we have:
\begin{itemize}
\item[(i)] $w_\rho \in D^{1,2}(\R^{N+1}_+;t^{1-2s})$, and $f_\rho \in L^{q_s}(\Omega_\rho)$.
\item[(ii)] If $w$ weakly solves the problem
\be\label{eq:smeq-1}
\left\{
\begin{aligned}
\div(t^{1-2s}\n w)&=0&&\quad \textrm{ in }\R^{N+1}_+,\\
- c_{N,s}\lim \limits_{t \to 0} t^{1-2s}\frac{\de w}{\de t} &= f &&  \quad \textrm{ on } \Omega,
\end{aligned}
\right.
\ee
then $w_\rho$ weakly solves the problem
\be\label{eq:smeq-1-rho}
\left\{
\begin{aligned}
\div(t^{1-2s}\n w_\rho)&=0&&\quad \textrm{ in }\R^{N+1}_+,\\
- c_{N,s}\lim \limits_{t \to 0} t^{1-2s}\frac{\de w_\rho}{\de t} &= f_\rho  && \quad \textrm{ on } \Omega_\rho.
\end{aligned}
\right.
\ee
\end{itemize}
\end{Proposition}

\proof
Let $w \in D^{1,2}(\R^{N+1}_+;t^{1-2s})$ and $f \in
L^{q_s}(\Omega)$. Note that $w_\rho(z)=\rho^{2s-N}w_1(\frac{z}{\rho^2})$
and $f_\rho(x)=\rho^{-(N+2s)}w_1(\frac{x}{\rho^2})$ for every $\rho>0$,
$z \in \R^{N+1}_+$ and $x \in \R^N \setminus \{0\}$. Hence it suffices to
prove the claims in the case $\rho=1$, and we put $\tilde w=w_1$,
$\tilde f= f_1$ and $\tilde \Omega= \Omega_1$. Recalling the properties of
the map $\kappa$ defined in (\ref{eq:15}), we then find
$$
\int_{\tilde \Omega} |\tilde f|^{q_s}\,dx = \int_{\tilde
  \Omega}|x|^{-2N}
|f(\frac{x}{|x|^2})|^{q_s}\,dx=\int_{\tilde
  \Omega}|J_\kappa| |f \circ \kappa|^{q_s}\,dx= \int_{\Omega}|f|^{q_s}\,dx.
$$
To simplify the  notations, we set
$$
\t: \overline{\R^{N+1}_+} \setminus \{0\} \to \overline{\R^{N+1}_+} \setminus \{0\},\qquad \t(z)= \frac{z}{|z|^2},
$$
so that the restriction of $\t$ to $\R^N \setminus \{0\}$ coincides
with $\kappa$. We note that the Jacobian $J_\t$ of $\t$ satisfies
$$
J_\t^T(z) J_\t(z)= |z|^{-4} I,
$$
where $I$ denotes the $(n+1) \times (n+1)$-identity matrix, and $\det
J_\t(z)= |z|^{-2N-2}$ for every $z \in \R^{N+1}_+$.\\
Next, we write $\tilde w= g \circ \t$ with $
g(z)=|z|^{N-2s}w(z)$. Moreover, we let
$\phi \in C^{\infty}_c(\overline {\R^{N+1}_+} \setminus \{0\})$ be
arbitrary, and define $\tilde \phi \in C^{\infty}_c(\overline
{\R^{N+1}_+} \setminus \{0\})$
by $\tilde \phi= h \circ \t$ with $h(z)=|z|^{N-2s}\phi(z)$. Considering
first the special case where $w \in
C^{\infty}_c(\overline {\R^{N+1}_+} \setminus \{0\})$, we then calculate
\begin{align*}
\int_{\R^{N+1}_+} t^{1-2s} \nabla \tilde w \nabla \tilde \phi\,dz &=\int_{\R^{N+1}_+} t^{1-2s} [J_\t(z) \nabla
g(\t(z)) ]\cdot
[J_\t(z)  \nabla h(\t(z))] \,dz \\
&=  \int_{\R^{N+1}_+}  t^{1-2s}|z|^{-4} \nabla g(\t(z)) \nabla  h(\t(z))\,dz \\
&=  \int_{\R^{N+1}_+} |z|^{-2N-2}  \Bigl(\frac{t}{|z|^2}\Bigr)^{1-2s}
|z|^{2(N-2s)} \nabla g(\t(z)) \nabla h(\t(z))\,dz \\
&= \int_{\R^{N+1}_+} |\det J_\t(z)| \Bigl(\frac{t}{|z|^2}\Bigr)^{1-2s}
|\t(z)|^{2(2s-N)} \nabla g(\t(z)) \nabla h(\t(z))\,dz \\
&= \int_{\R^{N+1}_+} t^{1-2s}
|z|^{2(2s-N)} \nabla g(z) \nabla h(z)\,dz.
\end{align*}
Noting that
$$
\nabla g(z)= (N-2s)|z|^{N-2s-2}z w(z)+ |z|^{N-2s}\nabla w(z)
$$
and
$$
\nabla h(z)= (N-2s)|z|^{N-2s-2}z \phi(z)+ |z|^{N-2s}\nabla \phi(z),
$$
we then conclude that
$$
\int_{\R^{N+1}_+} t^{1-2s} \nabla \tilde w \nabla \tilde \phi \,dz = \int_{\R^{N+1}_+} t^{1-2s} \nabla
w \nabla \phi\,dz  + I_1 + I_2 + I_3
$$
with
\begin{align*}
I_1&= (N-2s)^2 \int_{\R^{N+1}_+} t^{1-2s} |z|^{-2}w(z) \phi(z)\,dz,\\
I_2&= (N-2s) \int_{\R^{N+1}_+} t^{1-2s} |z|^{-2}w(z) z \nabla \phi(z)\,dz,\\
I_3&= (N-2s) \int_{\R^{N+1}_+} t^{1-2s} |z|^{-2}\phi(z) z \nabla w(z)\,dz.
\end{align*}
Since ${\rm div}_z [t^{1-2s}|z|^{-2}z] = (N-2s)t^{1-2s}|z|^{-2}$, it
follows that
$$
I_1+I_2+I_3 = (N-2s)\int_{\R^{N+1}_+} \div_z \Bigl(t^{1-2s}|z|^{-2}z w(z)
\phi(z)\Bigr)\,dz = 0
$$
and therefore
\begin{equation}
  \label{eq:6}
\int_{\R^{N+1}_+} t^{1-2s} \nabla \tilde w \nabla \tilde \phi\,dz = \int_{\R^{N+1}_+} t^{1-2s} \nabla
w \nabla \phi\,dz.
\end{equation}
By \cite{Froh}, we have  that
$C^{\infty}_c(\overline {\R^{N+1}_+} \setminus \{0\})$ is dense in
$D^{1,2}(\R^{N+1}_+;t^{1-2s})$ thus we deduce that (\ref{eq:6})
also holds for arbitrary $w, \phi \in D^{1,2}(\R^{N+1}_+;t^{1-2s})$,
while $\tilde w, \tilde \phi$ are also contained in
$D^{1,2}(\R^{N+1}_+;t^{1-2s})$. In particular, (i) is proved.\\
Moreover, (\ref{eq:6}) implies that $\tilde w$ is weakly
$L_s$-harmonic if $w$ is weakly $L_s$-harmonic. In addition, considering
the traces of $w$ and $\tilde w$ respectively,
Corollary~\ref{sec:some-preliminaries-2} implies that $(-\Delta)^s
\tilde w= \tilde f$ in distributional sense in $\tilde \Omega$ if
$(-\Delta)^s w= f$ in distributional sense in $\Omega$. Hence (ii)
follows from Lemma~\ref{sec:some-preliminaries-4}.
\QED

We will need the following version of a strong maximum
principle which is essentially a reformulation of \cite[Proposition 4.11]{CS}.

\begin{Lemma}\label{lem:smmp}
Let  $E$ be an open subset of $\R^N$, and let $w \in  D^{1,2}(\R^{N+1}_+;t^{1-2s})$ be a weak solution
of
$$
\left\{
\begin{aligned}
\div(t^{1-2s}\n w) &= 0&&\quad \textrm{in } \R^{N+1}_+,\\
- c_{N,s}\lim \limits_{t \to 0^+} t^{1-2s}\frac{\de w}{\de t}&= g&&  \quad \textrm{ on } E
\end{aligned}
\right.
$$
for some $g \in L^{q_s}(E) \cap C(E)$. Suppose furthermore that $w$ is
continuous and nonnegative on $E \times [0,r]$ for some $r>0$, and that
\begin{equation}
  \label{eq:17}
\text{$g(x) \ge 0$ for every $x \in E$ with $w(x)=0$.}
\end{equation}
If $w \not \equiv 0$ in $E$, then
$w$ is strictly positive in $E$ and
therefore $\inf \limits_K {w}>0$ for any compact set $K\subset E$.
\end{Lemma}

\proof
If $w \not \equiv 0$ on $E$, then $w>0$ in $E \times (0,r)$, since
$w$ is $L_s$-harmonic and nonnegative in this set. Suppose by
contradiction that $w(x_0)=0$ for some $x_0 \in E$. Then $g(x_0)<0$ by
\cite[Proposition 4.11]{CS}, which contradicts (\ref{eq:17}).
\QED

We will also need the following ''small volume'' maximum principle:

\begin{Lemma}\label{lem:smvmp}
Let $\gamma>0$. Then there exists $\delta=\delta(N,s,\gamma)>0$ with
the following property. If
\begin{itemize}
\item[(i)] $F \subset \R^{N+1}_+$ is an open subset with $\de F\cap\R^N
  \not= \varnothing$,
\item[(ii)] $E$ is a bounded open subset of $\R^N$ with $E \subset \partial F$,
\item[(iii)] $c\in L^\infty(E)$ is given with $\|c\|_{L^\infty(E)} \le \gamma$,
\item[(iv)] $w\in D^{1,2}(\R^{N+1}_+;t^{1-2s})$ is a weak solution of 
\be\label{eq:smeq}
\left\{
\begin{aligned}
\div(t^{1-2s}\n w) &\le 0\qquad &&\textrm{in }\R^{N+1}_+,\\
- c_{N,s} \lim \limits_{t \to 0^+} t^{1-2s}\frac{\de w}{\de t}& \ge c(x) w   && \qquad \textrm{on } E,
\end{aligned}
\right.
\ee
i.e.,
\begin{equation}
  \label{eq:9}
c_{N,s} \int_{F}t^{1-2s} \nabla w \nabla \phi\,dz \ge  \int_{E} c(x)w
\phi\,dx 
\end{equation}
for all nonnegative $\phi \in D(E,s)$,
\item[(v)] $w$ is continuous on $\overline F$ and satisfies $w \ge 0$ on
  $\partial F \setminus E$,
\item[(vi)] $|\{x \in E\::\:w<0\}| \le \d$,
\end{itemize}
then $w\geq0$ in $F$.
\end{Lemma}

\proof
We consider the function 
$$
v:\overline{\R^{N+1}_+} \to \R, \qquad  v(x)= \left
\{
\begin{aligned}
&\max(-w(x),0),&&\quad x \in \overline F,\\
&0&&\quad x \in \overline{\R^{N+1}_+} \setminus \overline F.
\end{aligned}
\right.
$$
It can be deduced from assumptions (i) and (ii) that the relative boundary
of $\overline F$ in $\overline{\R^{N+1}_+}$ is contained in
$\overline{\partial F \setminus E}$, so that $v$ is continuous on
$\overline {\R^{N+1}_+}$ by assumption (v). Moreover, $v \equiv 0$ on
  $\R^{N} \setminus E$. As a consequence, $v \in H^1_{loc}(\R^{N+1}_+)$ by
  \cite[Theorem 9.17 and Remark 19]{Br}, and
$$
\int_{\R^{N+1}_+}t^{1-2s}|\n v|^2\,dz \le
\int_{\R^{N+1}_+}t^{1-2s}|\n w|^2\,dz < \infty.
$$
Hence $v \in D^{1,2}(\R^{N+1}_+;t^{1-2s})$, and
Lemma~\ref{sec:some-preliminaries-3} implies that $v \in D(E,s)$. We
also note that combining (\ref{eq:12}) and (\ref{eq:13}) yields a
constant $C=C(N,s)>0$ such that
$$
\|v\|_{L^{2N/(N-2s)}(\R^N)}^2 \leq  C \int_{ \R^{N+1}_+}t^{1-2s}|\n v|^2dz.
$$
Applying (\ref{eq:9}) to $v$, we then obtain
\begin{align*}
c_{N,s} &\int_{\R^{N+1}_+}t^{1-2s}|\n v|^2 \,dz= - c_{N,s} \int_{\R^{N+1}_+}t^{1-2s}\n w
\cdot \n v\,
dz \le - \int_E c(x)w v dx\\
&= \int_E c(x)v^2 dx \le
\|c\|_{L^\infty(E)} \bigl|\{x \in
E\::\:w<0\}\bigr|^{N/2s}\|v\|_{L^{2N/(N-2s)}(\R^N)}^2\\
&\le \gamma \:\delta^{N/2s}\:C\: \int_{ \R^{N+1}_+}t^{1-2s}|\n v|^2dz.
\end{align*}
Hence, if $\delta < \Bigl( \frac{c_{N,s}}{\gamma C}\Bigr)^{2s/N}$, then
  $v \equiv 0$ in $\R^{N+1}_+$ and therefore $w \ge 0$ in $F$, as claimed.
\QED

\section{Proof of the main results}
\label{sec:proof-main-results}

\noindent
In this section we complete the proof of our main results. We begin
with the\\[0.2cm]
\noindent{\em Proof of Theorem \ref{th:neSupc}:}\\
We suppose by contradiction that there
 exists a nontrivial solution $u\in
C(\R^N \setminus\{0\}) \cap \calD^{s,2}(\O)$ of \eqref{eq:pblm}, and
we let $w \in D^{1,2}(\R^{N+1}_+;t^{1-2s})$ denote the corresponding $L_s$-harmonic extension of $u$
which weakly solves the problem
$$
\left \{
\begin{aligned}
\div (t^{1-2s} \nabla w)&=0&&\quad \text{in $\R^{N+1}_+$,}\\
-c_{N,s} \lim \limits_{t \to 0}t^{1-2s}w_t&=f(x,w) &&\quad \text{on $\Omega'$,}
\end{aligned}
\right.
$$
for every open subset $\O'\subset\O$ which is relatively compact in $\R^N \setminus \{0\}$.
Here, as before, we also write $w$ in place of $u$ for the trace on $\R^N$. We clearly have $w\in C(\ov{\R^{N+1}_+}\setminus\{0\})$.
Let  $R:= \sup\{|x|\::\:x \in \Omega\}>0$. For $\rho \in (0,R)$, we consider the Kelvin transform $w_\rho$ of $w$ as
defined in Proposition~\ref{sec:some-preliminaries-1}. We also put
$$
F_\rho:= \{z \in \R^{N+1}_+ \::\: |z| > \rho\},\quad E_\rho:= \{x \in
\Omega\::\: |x| > \rho\}\quad \text{and}\quad \tilde E_\rho:=
\{\frac{\rho^2 x}{|x|^2}\::\: x \in E_\rho\}.
$$
By definition of $R$ and since $\O$ is star-shaped with respect to the origin, $E_\rho$ and
$\tilde E_\rho$ are nonempty open subsets of $\Omega$ which are
relatively compact in $\R^N \setminus \{0\}$ for $\rho \in (0,R)$, so that the restrictions of the map $x \mapsto
f(x,w(x))$ to $E_\rho$ and $\tilde E_\rho$ are bounded and continuous.
By Proposition~\ref{sec:some-preliminaries-1}, the difference function
$v_\rho= w_\rho-w \in D^{1,2}(\R^{N+1}_+;t^{1-2s})$ weakly
solves the problem
$$
\left \{
\begin{aligned}
\div (t^{1-2s} \nabla v_\rho)&=0&&\quad \text{in $\R^{N+1}_+$,}\\
-c_{N,s} \lim \limits_{t \to 0}t^{1-2s}[v_\rho]_t&=g_\rho&&\quad \text{on $E_\rho$,}
\end{aligned}
\right.
$$
where $g_\rho$ is the bounded and continuous function on $E_\rho$
given by
$$
g_\rho(x)= \left(\frac{\rho}{|x|}\right)^{N+2s} f\left(\frac{\rho^2 x}{|x|^2},\left( \frac{\rho}{|x|}\right)^{2s-N}
 w_\rho(x) \right)- f(x,w(x)) .
$$
Moreover, by the supercriticality assumption \eqref{eq:supCri} we have
$$
g_\rho(x) \ge f(x,w_\rho(x))-f(x,w(x))= c_\rho(x) v_\rho(x) \qquad \text{for $x$ in $E_\rho$}
$$
with
$$
c_\rho: E_\rho \to \R, \qquad c_\rho(x)= \left\{
  \begin{aligned}
&\frac{f(x,w_\rho(x))-f(x, w(x))}{w_\rho(x)-w(x)}&&\quad \text{if $w(x)\neq w_\rho(x)$,}\\
&0&&\quad \text{if $w(x)= w_\rho(x)$.}
  \end{aligned}
\right.
$$
We also note that, since $f$ is assumed to be locally Lipschitz in its
second variable, we
have $c_\rho \in L^\infty(E_\rho)$ for $0<\rho<R$, Moreover, for $\tau \in (0,R)$ we have 
\begin{equation}
  \label{eq:8}
\gamma_\tau := \sup_{\rho \in [\tau,R)}\|c\|_{L^\infty(E_\rho)}\:<\:\infty 
\end{equation}
We now define
$$
\rho_*:= \inf \{\bar \rho \in (0,R)\::\: \text{$v_\rho \ge 0$ in $F_\rho$
  for $\rho \in [\bar \rho,R)$}\}.
$$
Since $|E_\rho \cap \{w_\rho<0\}|$ is small provided $\rho$ is
sufficiently close to $R$, Lemma~\ref{lem:smvmp} implies that
$\rho_*<R$. We claim that $\rho_*=0$. Indeed, suppose by contradiction
that $\rho_*>0$. By continuity, we then have $v_{\rho_*} \ge 0$ in
$F_{\rho_*}$. Moreover, $v_{\rho_*} \not \equiv 0$ in $E_{\rho_*}$ since
$$
v_{\rho_*}(x)>0 \qquad \text{for every $x \in \partial \Omega$
  with $|x|>\rho_*$.}
$$
By Lemma~\ref{lem:smmp}, we obtain $v_{\rho_*}>0$ in
$E_{\rho_*}$. We now fix $\tau \in (0,\rho_*)$ and choose $\delta>0$ as in
Lemma~\ref{lem:smvmp} according to $\gamma=\gamma_\tau$ as defined in
(\ref{eq:8}). Moreover, we choose a compact set $K \subset E_{\rho_*}$ such
that $|E_{\rho_*}\setminus K|<\delta$. Then $\inf \limits_{K} w_{\rho_*}>0$, and by continuity we also have
$$
K \subset E_\rho,\qquad |E_{\rho}\setminus K|<\delta \qquad
\text{and}\qquad \inf \limits_{K} w_{\rho}>0
$$
for $\rho \in (\tau,\rho_*)$ sufficiently close to $\rho_*$. Therefore
Lemma~\ref{lem:smvmp} implies that $v_\rho \ge 0$ in $F_\rho$
  for $\rho \in (\tau,\rho_*)$ sufficiently close to $\rho_*$. This contradicts 
the definition of $\rho_*$. We conclude that
$\rho_*=0$, as claimed. As a consequence, for every $x\in \O$ and
$x\neq 0$ we have
\be\label{eq:momrho}
\left(\frac{\rho}{|x|}\right)^{N-2s}w\left(\frac{\rho^2
    x}{|x|^2}\right)\geq w(x) \quad \text{for all $\rho\in(0,|x|)$.}
\ee
Furthermore, since $w \in \calD^{s,2}(\Omega) \subset
L^{\frac{2N}{N-2s}}(\R^N)$ we have
$$
\int_{ S^{N-1}}\int_0^\infty w^{\frac{2N}{N-2s}}(r\s)drd\s = \int_{\R^N}w^{\frac{2N}{N-2s}}\,dx <\infty
$$
and therefore, by Fubini's theorem,
\be\label{eq:Sob-sf}
 \int_0^\infty w^{\frac{2N}{N-2s}}(r\s_0)dr<\infty \qquad \text{for
   a.e. $\s_0 \in S^{N-1}$.}
\ee
We now pick $\s_0 \in S^{N-1}$ and $r_0>0$ such that $r_0\s_0\in \O$ and \eqref{eq:Sob-sf}
holds for $\s_0$. By \eqref {eq:momrho} we then have
$$
\left(\frac{\rho}{r_0}\right)^{N-2s}w\left(\frac{\rho^2
    \s_0}{r_0}\right)\geq w(r_0\s_0 )>0 \qquad \text{for $\rho \in (0,r_0)$}
$$
and consequently
$$
w(r\s_0)\geq C r^{(2s-N)/2}\qquad \text{for $r\in(0,r_0)$ with a constant $C>0$.}
$$
 This implies
$$
\int_0^{r_0} w^{\frac{2N}{N-2s}}(r\s_0)dr=\infty
$$
contrary to \eqref{eq:Sob-sf}. The contradiction shows that there
does not exist a nontrivial solution $u\in
C(\R^N \setminus\{0\}) \cap \calD^{s,2}(\O)$ of \eqref{eq:pblm}
under the assumptions of Theorem~\ref{th:neSupc}, as claimed.
\QED

\noindent{\em Proof of Corollary~\ref{cor-poho}:}\\
Problem \eqref{eq:potV} is a special case of \eqref{eq:pblm} with
$f(x,u)= u^p-V(x)u$, and for this nonlinearity we calculate
$$
H_f(x,u)= \Bigl(p-\frac{N+2s}{N-2s}\Bigr)u^p +
\frac{4u}{N-2s}\Bigl(sV(x)+\frac{1}{2} x \cdot \nabla V(x)\Bigr)
$$
so that (\ref{eq:3}) is satisfied by the assumptions on $p$ and $V$.
Moreover, any nontrivial,
nonnegative solution of \eqref{eq:potV} is strictly positive in
$\Omega \setminus \{0\}$, which follows by applying
Lemma~\ref{lem:smmp} to the $L_s$-harmonic extension of $u$ and the
sets $E_\eps:= \{x \in \Omega \::\: |x|>\eps\}$ for $\eps>0$
small. Hence nontrivial,
nonnegative solutions of \eqref{eq:potV} do not exist by Theorem~\ref{th:neSupc}.
\QED

\noindent{\em Proof of Corollary~\ref{cor:hardy}:}\\
Problem \eqref{eq:pdHard} is a special case of \eqref{eq:potV} with
$V(x)=\gamma |x|^{-2s}$, so the result follows from Corollary~\ref{cor-poho}.
\QED

\noindent{\em Proof of Corollary~\ref{cor:pdSW}:}\\
Problem \eqref{eq:pdSW} is a special case of \eqref{eq:pblm} with
$f(x,u)= |x|^{-\sigma} u^p$, and for this nonlinearity we calculate
$$
H_f(x,u)= \Bigl(p-\frac{N+2s-2\sigma}{N-2s}\Bigr)|x|^{-\sigma}u^p.
$$
Hence (\ref{eq:3}) is satisfied by the assumptions on $p$ and $\sigma$.
Moreover, by the same argument as in the proof of
Corollary~\ref{cor-poho} above, nontrivial and  nonnegative solutions
of \eqref{eq:pdSW} must be strictly positive in
$\Omega \setminus \{0\}$ and therefore can not exist by Theorem~\ref{th:neSupc}.    \QED

\noindent We finally give the proof of our nonexistence result in
(unbounded) domains being star-shaped at infinity.\\[0.2cm]
\noindent \textit{Proof of Theorem \ref{th:ne-stInf}:}\\
The definition of star-shapedness at infinity implies that, after a
suitable translation, the image $\tilde \Omega:= \kappa(\Omega)$ of
the domain  $\Omega$ under the map $\kappa$ defined in (\ref{eq:15}) is
star-shaped with respect to $0 \in \partial \tilde
\Omega$. Moreover, if $u \in \calD^{s,2}(\tilde \Omega) \cap C(\R^N)$ is a nonnegative solution of \eqref{eq:ustf-stU}, then
Corollary~\ref{sec:some-preliminaries-2} implies that $\tilde u = K u \in \calD^{s,2}(\tilde \Omega) \cap C(\R^N
\setminus \{0\})$ solves $(-\Delta)^s \tilde u = |x|^{-\sigma}|\tilde u|^p$
in $\tilde \Omega$ with $\sigma= N+2s-p(N-2s)$. Since the assumption
$p \le
\frac{N+2s}{N-2s}$ yields $p \ge \frac{N+2s-2\sigma}{N-2s}$,
Corollary~\ref{cor:pdSW} implies that $\tilde u \equiv 0$ and  hence also
$u \equiv 0$, as claimed.
\QED


    \label{References}

   \end{document}